\documentclass{amsart}[12pt]
\begin{document}
\begin{center}
\textbf{\Large{LINEAR INDEPENDENT SOLUTIONS AND OPERATIONAL REPRESENTATIONS FOR HYPERGEOMETRIC FUNCTIONS OF FOUR VARIABLES}}
\medskip\\
\textbf{Maged G. Bin-Saad}$^1$ {\bf and Anvar Hasanov}$^2$\\
\medskip
$^1$Department of Mathematics, Aden University, Khormaksar, P.O. Box 6014, Aden, Yemen\\
$^2$Institute of Mathematics, 29 F. Hodjaev Street, Tashkent 700125, Uzbekistan \\
E-Mails: mgbinsaad@yahoo.com, anvarhasanov@yahoo.com
\end{center}

\medskip

{\bf Abstract}

  In investigation of boundary-value problems for certain partial differential equations arising in applied mathematics, we often need to study the solution of  system   of partial differential equations satisfied by hypergeometric functions and  find explicit linearly independent solutions for the system.  Here  we choose the Exton function $K_{2}$ among his 21 functions
to show how to find the linearly independent solutions of partial differential equations satisfied by this function $K_{2}$. Based upon the classical derivative and integral operators we introduce a new operational images for hypergeometric function $K_{2}$.
By means of these operational images a number of finite series and decomposition formulas are then fund.

\vspace{5 mm }
 \vskip 5 mm
\section*{\bf 1. Introduction and Preliminaries}
Solutions of many applied problems involving thermal conductivity and dynamics, electromagnetic oscillation and aerodynamics, quantum mechanics and potential theory are obtainable with the help of hypergeometric (higher and special or transcendent) functions \cite{Bers, Fran, Niuk,Loho}.
 Functions of such kind are often referred to as special functions of mathematical physics. They mainly  appear in the solution of partial differential equations which are dealt with  harmonic analysis method (see \cite{Er-Ma-Ob-Tr-1}).
In view of various applications, it is interesting in itself and seems to be very important to conduct a continuous research of multiple hypergeometric functions. For instance, in \cite{Sr-Ka-85}, a comprehensive list of hypergeometric functions of three variables as many as 205 is recorded, together with their regions of convergence.
	It is noted that Riemann's functions and the fundamental solutions of the degenerate second-order partial differential equations are expressible by means of hypergeometric functions of several variables (see \cite{Alti,Ba-Ge-99,Ba-Ge-02,Ba-Ge-05,Frya,Hasa-07-a,Hasa-07-b,Hasa-07-c, Hasa-08, Ha-Ra-Tu,Ha-Ka-09,Mcco, Sa-Ha-08}). Therefore, in investigation of boundary-value problems for these partial differential equations, we need to study the solution of the system of hypergeometric functions and  find explicit linearly independent solutions (see \cite{Hasa-07-a,Hasa-07-b,Hasa-07-c, Hasa-08, Ha-Ra-Tu,Ha-Ka-09,Sa-Ha-08}).
Exton (\cite{Ex-76}, pp. 78-79) introduced 21 complete hypergeometric functions ${{K}_{1}},{{K}_{2}},...,{{K}_{21}}$ of four variables. In \cite{Sh-Pa-89} Sharma and Parihar introduced 83 complete hypergeometric functions $F_{1}^{\left( 4 \right)},\,F_{2}^{\left( 4 \right)},\,...,F_{83}^{\left( 4 \right)}\,$of four variables. It is remarkable that out of these 83 functions, the following 19 functions had already appeared in the work of Exton \cite{Ex-76} in the different notations:
$$
  F_{9}^{\left( 4 \right)}={{K}_{1}},\,\,F_{1}^{\left( 4 \right)}={{K}_{2}},\,\,F_{38}^{\left( 4 \right)}={{K}_{3}},\,\,F_{10}^{\left( 4 \right)}={{K}_{4}},\,\,F_{2}^{\left( 4 \right)}={{K}_{5}},\,\,F_{59}^{\left( 4 \right)}={{K}_{6}},\,\, \\
$$
$$ F_{39}^{\left( 4 \right)}={{K}_{7}},\,\,F_{11}^{\left( 4 \right)}={{K}_{8}},\,\,F_{12}^{\left( 4 \right)}={{K}_{9}},\,\,F_{3}^{\left( 4 \right)}={{K}_{10}},\,\,F_{60}^{\left( 4 \right)}={{K}_{11}},\,\,F_{40}^{\left( 4 \right)}={{K}_{12}},\,\, \\
$$
$$ F_{13}^{\left( 4 \right)}={{K}_{13}},\,\,F_{77}^{\left( 4 \right)}={{K}_{14}},\,\,F_{78}^{\left( 4 \right)}={{K}_{15}},\,\,F_{79}^{\left( 4 \right)}={{K}_{16}},\,\,F_{82}^{\left( 4 \right)}={{K}_{19}},\,\,F_{81}^{\left( 4 \right)}={{K}_{20}},\,\,F_{83}^{\left( 4 \right)}={{K}_{21}}.\,\,
$$
Each quadruple hypergeometric function is of the form
$${{F}^{\left( 4 \right)}}\left( \cdot  \right)=\sum\limits_{m,n,p,q=0}^{\infty }{{}}\Delta \left( m,n,p,q \right)\frac{{{x}^{m}}}{m!}\frac{{{y}^{n}}}{n!}\frac{{{z}^{p}}}{p!}\frac{{{t}^{q}}}{q!}
$$
where $\Delta \left( m,n,p,q \right)$ is a certain sequence of complex parameters and there are twelve parameters in each function ${{F}^{\left( 4 \right)}}\left( \cdot  \right)$.
Here, for an example, we choose the Exton function $K_{2}$ among his twenty one functions
$$
K_{2}\left( a,b,c;{{e}_{1}},{{e}_{2}},{{e}_{3}},{{e}_{4}};x,y,z,t \right)=\sum\limits_{m,n,p,q=0}^{\infty }{{}}\frac{{{\left( a \right)}_{m+n+p+q}}{{\left( b \right)}_{m+n+q}}{{\left( c \right)}_{p}}}{{{\left( {{e}_{1}} \right)}_{m}}{{\left( {{e}_{2}} \right)}_{n}}{{\left( {{e}_{3}} \right)}_{p}}{{\left( {{e}_{4}} \right)}_{q}}}\frac{{{x}^{m}}}{m!}\frac{{{y}^{n}}}{n!}\frac{{{z}^{p}}}{p!}\frac{{{t}^{q}}}{q!}\eqno (1.1)
$$
to find the linearly independent solutions of partial differential equations satisfied by this function $K_{2}$.

\section*{2. The system of partial differential equations for  $K_{2}$}

According to the theory of multiple hypergeometric functions (see \cite{Ap-Ka}),  the system of partial differential equations for the Sharma and Parihar hypergeometric function $K_{2}$ is readily seen to be given as follows:
$$
\left\{ \begin{matrix}
   \left( {{e}_{1}}+x\frac{\partial }{\partial x} \right)\left( x\frac{\partial }{\partial x}+1 \right){{x}^{-1}}u-\left( a+x\frac{\partial }{\partial x}+y\frac{\partial }{\partial y}+z\frac{\partial }{\partial z}+t\frac{\partial }{\partial t} \right)\left( b+x\frac{\partial }{\partial x}+y\frac{\partial }{\partial y}+t\frac{\partial }{\partial t} \right)u=0  \\
   \left( {{e}_{2}}+y\frac{\partial }{\partial y} \right)\left( y\frac{\partial }{\partial y}+1 \right){{y}^{-1}}u-\left( a+x\frac{\partial }{\partial x}+y\frac{\partial }{\partial y}+z\frac{\partial }{\partial z}+t\frac{\partial }{\partial t} \right)\left( b+x\frac{\partial }{\partial x}+y\frac{\partial }{\partial y}+t\frac{\partial }{\partial t} \right)u=0  \\
   \,\,\,\,\,\,\,\,\,\,\,\,\,\,\,\,\,\,\,\,\,\,\,\,\,\,\,\,\,\,\,\,\,\,\,\,\,\,\, \left( {{e}_{3}}+z\frac{\partial }{\partial z} \right)\left( z\frac{\partial }{\partial z}+1 \right){{z}^{-1}}u-\left( a+x\frac{\partial }{\partial x}+y\frac{\partial }{\partial y}+z\frac{\partial }{\partial z}+t\frac{\partial }{\partial t} \right)\left( c+z\frac{\partial }{\partial z} \right)u=0  \\
   \,\,\,\,\,\,\,\,\, \left( {{e}_{4}}+t\frac{\partial }{\partial t} \right)\left( t\frac{\partial }{\partial t}+1 \right){{t}^{-1}}u-\left( a+x\frac{\partial }{\partial x}+y\frac{\partial }{\partial y}+z\frac{\partial }{\partial z}+t\frac{\partial }{\partial t} \right)\left( b+x\frac{\partial }{\partial x}+y\frac{\partial }{\partial y}+t\frac{\partial }{\partial t} \right)u=0  \\
\end{matrix} \right.\eqno (2.1)
$$
where $u=K_{2}\left( a,b,c;{{e}_{1}},{{e}_{2}},{{e}_{3}},{{e}_{4}};x,y,z,t \right)$. \\
Now by making use of some elementary calculations,  we find the following system of second order partial differential equations:
$$
\left\{ \begin{matrix}
    x\left( 1-x \right){{u}_{xx}}-{{y}^{2}}{{u}_{yy}}-{{t}^{2}}{{u}_{tt}}-2xy{{u}_{xy}}-xz{{u}_{xz}}-2xt{{u}_{xt}}-yz{{u}_{yz}}-2yt{{u}_{yt}}-zt{{u}_{zt}} \\
 \,\,\,\,\,\,\,\,\,\,\,\,\,\,\,\,\,\,\,\,\,\,\,\,\,\,\,\,\,\,\,\,\,\,\,\,\,\,\,\,\,\,\, +\left[ {{e}_{1}}-\left( a+b+1 \right)x \right]{{u}_{x}}-\left( a+b+1 \right)y{{u}_{y}}-bz{{u}_{z}}-\left( a+b+1 \right)t{{u}_{t}}-abu=0 \\
     y\left( 1-y \right){{u}_{yy}}-{{x}^{2}}{{u}_{xx}}-{{t}^{2}}{{u}_{tt}}-2xy{{u}_{xy}}-xz{{u}_{xz}}-2xt{{u}_{xt}}-yz{{u}_{yz}}-2yt{{u}_{yt}}-zt{{u}_{zt}} \\
  \,\,\,\,\,\,\,\,\,\,\,\,\,\,\,\,\,\,\,\,\,\,\,\,\,\,\,\,\,\,\,\,\,\,\,\,\,\,\,\,\,\,\,\, -\left( a+b+1 \right)x{{u}_{x}}+\left[ {{e}_{2}}-\left( a+b+1 \right)y \right]{{u}_{y}}-bz{{u}_{z}}-\left( a+b+1 \right)t{{u}_{t}}-abu=0 \\
     z\left( 1-z \right){{u}_{zz}}-xz{{u}_{xz}}-yz{{u}_{yz}}-zt{{u}_{zt}}-cx{{u}_{x}}-cy{{u}_{y}}+\left[ {{e}_{3}}-\left( a+c+1 \right)z \right]{{u}_{z}}-ct{{u}_{t}}-acu=0  \\
     t\left( 1-t \right){{u}_{tt}}-{{x}^{2}}{{u}_{xx}}-{{y}^{2}}{{u}_{yy}}-2xy{{u}_{xy}}-xz{{u}_{xz}}-2xt{{u}_{xt}}-yz{{u}_{yz}}-2yt{{u}_{yt}}-zt{{u}_{zt}} \\
 \,\,\,\,\,\,\,\,\,\,\,\,\,\,\,\,\,\,\,\,\,\,\,\,\,\,\,\,\,\,\,\,\,\,\,\,\,\,\,\,\,\,\,\,\, -\left( a+b+1 \right)x{{u}_{x}}-\left( a+b+1 \right)y{{u}_{y}}-bz{{u}_{z}}+\left[ {{e}_{4}}-\left( a+b+1 \right)t \right]{{u}_{t}}-abu=0 \\
\end{matrix} \right.\eqno (2.2)
$$
It is noted that three equations of the system (2.2) are simultaneous, because the hypergeometric function $K_{2}$
satisfies the system. Now in order to find the linearly independent solutions of the system (2.2) we consider $u$ as in the form $u={{x}^{\alpha }}{{y}^{\beta }}{{z}^{\gamma }}{{t}^{\delta }}w$ is an unknown function, and $\alpha$, $\beta$, $\gamma$
and $\delta $ are constants which are to be determined. So, substituting $u={{x}^{\alpha }}{{y}^{\beta }}{{z}^{\gamma }}{{t}^{\delta }}w$
into the system (2.2),  we obtain
\vskip 0.1 mm
$$
\left\{ \begin{matrix}
    x\left( 1-x \right){{w}_{xx}}-{{y}^{2}}{{w}_{yy}}-{{t}^{2}}{{w}_{tt}}-2xy{{w}_{xy}}-xz{{w}_{xz}}-2xt{{w}_{xt}}-yz{{w}_{yz}}-2yt{{w}_{yt}}-zt{{w}_{zt}} \\
  +\left[ {{E}_{1}}-\left( A+B+1 \right)x \right]{{w}_{x}}-\left( A+B+1 \right)y{{w}_{y}}-Bz{{w}_{z}}-\left( A+B+1 \right)t{{w}_{t}} \\
  \,\,\,\,\,\,\,\,\,\,\,\,\,\,\,\,\,\,\,\,\,\,\,\,\,\,\,\,\,\,\,\,\,\,\,\,\,\,\,\,\,\,\,\,\,\,\,\,\,\,\,\,\,\,\,\,\,\,\,\,\,\,\,\,\,\,\,
  \,\,\,\,\,\,\,\,\,\,\,\,\,\,\,\,\,\,\,\,\,\,\,\,\,\,\,\,\,\,\,\,\,\,\,\,\,\,\,\,\,\,\,\,\,\,\,\,\,\,\,\,\,\,\,\,\,\,\,\,\,\,\,\,\,\,\,\,\,\,\,\,\,
  \,\,\,\,\,\,\,\, +\left[ \alpha \left( \alpha -1+{{e}_{1}} \right){{x}^{-1}}-AB \right]w=0 \\
       y\left( 1-y \right){{w}_{yy}}-{{x}^{2}}{{w}_{xx}}-{{t}^{2}}{{w}_{tt}}-2xy{{w}_{xy}}-xz{{w}_{xz}}-2xt{{w}_{xt}}-yz{{w}_{yz}}-2yt{{w}_{yt}}-zt{{w}_{zt}} \\
  -\left( A+B+1 \right)x{{w}_{x}}+\left[ {{E}_{2}}-\left( A+B+1 \right)y \right]{{w}_{y}}-Bz{{w}_{z}}-\left( A+B+1 \right)t{{w}_{t}} \\
  \,\,\,\,\,\,\,\,\,\,\,\,\,\,\,\,\,\,\,\,\,\,\,\,\,\,\,\,\,\,\,\,\,\,\,\,\,\,\,\,\,\,\,\,\,\,\,\,\,\,\,\,\,\,\,\,\,\,\,\,\,\,\,\,\,\,\,
  \,\,\,\,\,\,\,\,\,\,\,\,\,\,\,\,\,\,\,\,\,\,\,\,\,\,\,\,\,\,\,\,\,\,\,\,\,\,\,\,\,\,\,\,\,\,\,\,\,\,\,\,\,\,\,\,\,\,\,\,\,\,\,\,\,\,\,\,\,\,\,\,\,
  \,\,\,\,\,\,\,\,\,\,\,\, +\left[ \beta \left( \beta -1+{{e}_{2}} \right){{y}^{-1}}-AB \right]w=0 \\
        z\left( 1-z \right){{w}_{zz}}-xz{{w}_{xz}}-yz{{w}_{yz}}-zt{{w}_{zt}}-Cx{{w}_{x}}-Cy{{w}_{y}}+\left[ {{E}_{3}}-\left( A+C+1 \right)z \right]{{w}_{z}}-Ct{{w}_{t}} \\
  \,\,\,\,\,\,\,\,\,\,\,\,\,\,\,\,\,\,\,\,\,\,\,\,\,\,\,\,\,\,\,\,\,\,\,\,\,\,\,\,\,\,\,\,\,\,\,\,\,\,\,\,\,\,\,\,\,\,\,\,\,\,\,\,\,\,\,
  \,\,\,\,\,\,\,\,\,\,\,\,\,\,\,\,\,\,\,\,\,\,\,\,\,\,\,\,\,\,\,\,\,\,\,\,\,\,\,\,\,\,\,\,\,\,\,\,\,\,\,\,\,\,\,\,\,\,\,\,\,\,\,\,\,\,\,\,\,
  \,\,\,\,\,\,\,\,\,\,\,\,\,\,\,\,\,\,\,\, +\left[ \gamma \left( \gamma -1+{{e}_{3}} \right){{z}^{-1}}-AC \right]w=0 \\
        t\left( 1-t \right){{w}_{tt}}-{{x}^{2}}{{w}_{xx}}-{{y}^{2}}{{w}_{yy}}-2xy{{w}_{xy}}-xz{{w}_{xz}}-2xt{{w}_{xt}}-yz{{w}_{yz}}-2yt{{w}_{yt}}-zt{{w}_{zt}} \\
  -\left( A+B+1 \right)x{{w}_{x}}-\left( A+B+1 \right)y{{w}_{y}}-Bz{{w}_{z}}+\left[ {{E}_{4}}-\left( A+B+1 \right)t \right]{{w}_{t}}\,\, \\
  \,\,\,\,\,\,\,\,\,\,\,\,\,\,\,\,\,\,\,\,\,\,\,\,\,\,\,\,\,\,\,\,\,\,\,\,\,\,\,\,\,\,\,\,\,\,\,\,\,\,\,\,\,\,\,\,\,\,\,\,\,\,\,\,\,\,\,\,\,
  \,\,\,\,\,\,\,\,\,\,\,\,\,\,\,\,\,\,\,\,\,\,\,\,\,\,\,\,\,\,\,\,\,\,\,\,\,\,\,\,\,\,\,\,\,\,\,\,\,\,\,\,\,\,\,\,\,\,\,\,\,\,\,\,\,\,\,
   \,\,\,\,\,\,\,\,\,\,\,\,\,\,\,\,\,\,\, +\left[ \delta \left( \delta -1+{{e}_{4}} \right){{t}^{-1}}-AB \right]w=0 \\
  \end{matrix} \right.\eqno (2.3)
$$
where
$$
   A=\alpha +\beta +\gamma +\delta +a,\,\,B=\alpha +\beta +\delta +b,\,\,C=\gamma +c,\,\, {{E}_{1}}=2\alpha +{{e}_{1}},\,\,{{E}_{2}}=2\beta +{{e}_{2}},\,\,{{E}_{3}}=2\gamma +{{e}_{3}},\,\,{{E}_{4}}=2\delta +{{e}_{4}}.
$$
It is noted  that the system (2.3) is analogical to the system (2.2). Therefore, it is required that the following conditions should be satisfied:
 \vskip 1 mm
$$
\left\{ \begin{matrix}
   \alpha \left( \alpha -1+{{e}_{1}} \right)=0  \\
   \beta \left( \beta -1+{{e}_{2}} \right)=0  \\
   \gamma \left( \gamma -1+{{e}_{3}} \right)=0  \\
   \delta \left( \delta -1+{{e}_{4}} \right)=0  \\
\end{matrix} \right.\eqno (2.4)
$$
It is not difficult to see that the system (2.4) satisfies the following solutions:
$$
\begin{matrix}
   {} & 1 & 2 & 3 & 4 & 5 & 6 & 7 & 8  \\
   \alpha : & 0 & 1-{{e}_{1}} & 0 & 0 & 0 & 1-{{e}_{1}} & 1-{{e}_{1}} & 1-{{e}_{1}}  \\
   \beta : & 0 & 0 & 1-{{e}_{2}} & 0 & 0 & 1-{{e}_{2}} & 0 & 0  \\
   \gamma : & 0 & 0 & 0 & 1-{{e}_{3}} & 0 & 0 & 1-{{e}_{3}} & 0  \\
   \delta : & 0 & 0 & 0 & 0 & 1-{{e}_{4}} & 0 & 0 & 1-{{e}_{4}}  \\
\end{matrix}
$$
$$
\begin{matrix}
   {} & 9 & 10 & 11 & 12 & 13 & 14 & 15 & 16  \\
   \alpha : & 0 & 0 & 0 & 1-{{e}_{1}} & 1-{{e}_{1}} & 1-{{e}_{1}} & 0 & 1-{{e}_{1}}  \\
   \beta : & 1-{{e}_{2}} & 1-{{e}_{2}} & 0 & 1-{{e}_{2}} & 1-{{e}_{2}} & 0 & 1-{{e}_{2}} & 1-{{e}_{2}}  \\
   \gamma : & 1-{{e}_{3}} & 0 & 1-{{e}_{3}} & 1-{{e}_{3}} & 0 & 1-{{e}_{3}} & 1-{{e}_{3}} & 1-{{e}_{3}}  \\
   \delta : & 0 & 1-{{e}_{4}} & 1-{{e}_{4}} & 0 & 1-{{e}_{4}} & 1-{{e}_{4}} & 1-{{e}_{4}} & 1-{{e}_{4}}  \\
\end{matrix}\eqno (2.5)
$$
Finally, substituting all eight solutions (2.5) into (2.3),  we find the following linearly independent solutions of the system (2.2):
$$
{{u}_{1}}\left( x,y,z,t \right)=K_{2}\left( a,b,c;{{e}_{1}},{{e}_{2}},{{e}_{3}},{{e}_{4}};x,y,z,t \right),\eqno (2.6)
$$
$$
{{u}_{2}}\left( x,y,z,t \right)={{x}^{1-{{e}_{1}}}}K_{2}\left( 1-{{e}_{1}}+a,1-{{e}_{1}}+b,c;2-{{e}_{1}},{{e}_{2}},{{e}_{3}},{{e}_{4}};x,y,z,t \right),\eqno (2.7)
$$
$$
{{u}_{3}}\left( x,y,z,t \right)={{y}^{1-{{e}_{2}}}}K_{2}\left( 1-{{e}_{2}}+a,1-{{e}_{2}}+b,c;{{e}_{1}},2-{{e}_{2}},{{e}_{3}},{{e}_{4}};x,y,z,t \right),\eqno (2.8)
$$
$$
{{u}_{4}}\left( x,y,z,t \right)={{z}^{1-{{e}_{3}}}}K_{2}\left( 1-{{e}_{3}}+a,b,1-{{e}_{3}}+c;{{e}_{1}},{{e}_{2}},2-{{e}_{3}},{{e}_{4}};x,y,z,t \right),\eqno (2.9)
$$
$$
{{u}_{5}}\left( x,y,z,t \right)={{t}^{1-{{e}_{4}}}}K_{2}\left( 1-{{e}_{4}}+a,1-{{e}_{4}}+b,c;{{e}_{1}},{{e}_{2}},{{e}_{3}},2-{{e}_{4}};x,y,z,t \right),\eqno (2.10)
$$
$$
{{u}_{6}}\left( x,y,z,t \right)={{x}^{1-{{e}_{1}}}}{{y}^{1-{{e}_{2}}}}K_{2}\left( 2-{{e}_{1}}-{{e}_{2}}+a,2-{{e}_{1}}-{{e}_{2}}+b,c;2-{{e}_{1}},2-{{e}_{2}},{{e}_{3}},{{e}_{4}};x,y,z,t \right),\eqno (2.11)
$$
$$
{{u}_{7}}\left( x,y,z,t \right)={{x}^{1-{{e}_{1}}}}{{z}^{1-{{e}_{3}}}}K_{2}\left( 2-{{e}_{1}}-{{e}_{3}}+a,1-{{e}_{1}}+b,1-{{e}_{3}}+c;2-{{e}_{1}},{{e}_{2}},2-{{e}_{3}},{{e}_{4}};x,y,z,t \right),\eqno (2.12)
$$
$$
{{u}_{8}}\left( x,y,z,t \right)={{x}^{1-{{e}_{1}}}}{{t}^{1-{{e}_{4}}}}K_{2}\left( 2-{{e}_{1}}-{{e}_{4}}+a,2-{{e}_{1}}-{{e}_{4}}+b,c;2-{{e}_{1}},{{e}_{2}},{{e}_{3}},2-{{e}_{4}};x,y,z,t \right),\eqno (2.13)
$$
$$
{{u}_{9}}\left( x,y,z,t \right)={{y}^{1-{{e}_{2}}}}{{z}^{1-{{e}_{3}}}}K_{2}\left( 2-{{e}_{2}}-{{e}_{3}}+a,1-{{e}_{2}}+b,1-{{e}_{3}}+c;{{e}_{1}},2-{{e}_{2}},2-{{e}_{3}},{{e}_{4}};x,y,z,t \right),\eqno (2.14)
$$
$$
{{u}_{10}}\left( x,y,z,t \right)={{y}^{1-{{e}_{2}}}}{{t}^{1-{{e}_{4}}}}K_{2}\left( 2-{{e}_{2}}-{{e}_{4}}+a,2-{{e}_{2}}-{{e}_{4}}+b,c;{{e}_{1}},2-{{e}_{2}},{{e}_{3}},2-{{e}_{4}};x,y,z,t \right),\eqno (2.15)
$$
$$
{{u}_{11}}\left( x,y,z,t \right)={{z}^{1-{{e}_{3}}}}{{t}^{1-{{e}_{4}}}}K_{2}\left( 2-{{e}_{3}}-{{e}_{4}}+a,1-{{e}_{4}}+b,1-{{e}_{3}}+c;{{e}_{1}},{{e}_{2}},2-{{e}_{3}},2-{{e}_{4}};x,y,z,t \right),\eqno (2.16)
$$
$$
\begin{gathered}
  {{u}_{12}}\left( x,y,z,t \right)={{x}^{1-{{e}_{1}}}}{{y}^{1-{{e}_{2}}}}{{z}^{1-{{e}_{3}}}} \\
  \times K_{2}\left( 3-{{e}_{1}}-{{e}_{2}}-{{e}_{3}}+a,2-{{e}_{1}}-{{e}_{2}}+b,1-{{e}_{3}}+c;2-{{e}_{1}},2-{{e}_{2}},2-{{e}_{3}},{{e}_{4}};x,y,z,t \right),
\end{gathered}\eqno (2.17)
$$
$$
\begin{gathered}
  {{u}_{13}}\left( x,y,z,t \right)={{x}^{1-{{e}_{1}}}}{{y}^{1-{{e}_{2}}}}{{t}^{1-{{e}_{4}}}} \\
 \times K_{2}\left( 3-{{e}_{1}}-{{e}_{2}}-{{e}_{4}}+a,3-{{e}_{1}}-{{e}_{2}}-{{e}_{4}}+b,c;2-{{e}_{1}},2-{{e}_{2}},{{e}_{3}},2-{{e}_{4}};x,y,z,t \right),
\end{gathered}\eqno (2.18)
$$
$$
\begin{gathered}
   {{u}_{14}}\left( x,y,z,t \right)={{x}^{1-{{e}_{1}}}}{{z}^{1-{{e}_{3}}}}{{t}^{1-{{e}_{4}}}} \\
 \times K_{2}\left( 3-{{e}_{1}}-{{e}_{3}}-{{e}_{4}}+a,2-{{e}_{1}}-{{e}_{4}}+b,1-{{e}_{3}}+c;2-{{e}_{1}},{{e}_{2}},2-{{e}_{3}},2-{{e}_{4}};x,y,z,t \right),
\end{gathered}\eqno (2.19)
$$
$$
\begin{gathered}
  {{u}_{15}}\left( x,y,z,t \right)={{y}^{1-{{e}_{2}}}}{{z}^{1-{{e}_{3}}}}{{t}^{1-{{e}_{4}}}} \\
 \times K_{2}\left( 3-{{e}_{2}}-{{e}_{3}}-{{e}_{4}}+a,2-{{e}_{2}}-{{e}_{4}}+b,1-{{e}_{3}}+c;{{e}_{1}},2-{{e}_{2}},2-{{e}_{3}},2-{{e}_{4}};x,y,z,t \right),
\end{gathered}\eqno (2.20)
$$
$$
\begin{gathered}
 {{u}_{16}}\left( x,y,z,t \right)={{x}^{1-{{e}_{1}}}}{{y}^{1-{{e}_{2}}}}{{z}^{1-{{e}_{3}}}}{{t}^{1-{{e}_{4}}}} \\
 \times K_{2}\left( 4-{{e}_{1}}-{{e}_{2}}-{{e}_{3}}-{{e}_{4}}+a,3-{{e}_{1}}-{{e}_{2}}-{{e}_{4}}+b,1-{{e}_{3}}+c;2-{{e}_{1}},2-{{e}_{2}},2-{{e}_{3}},2-{{e}_{4}};x,y,z,t \right).
\end{gathered}\eqno (2.21)
$$
Moreover, it is seen that the global solution of the system  (2.2) is combined to be in the form
 $$
 u = \sum_{j=1}^{16}\, k_j \,  u_j \eqno (2.22)
 $$
where $ k_i$ $(i = 1,\, 2,\ldots,16)$ are constants.
\section*{3. Operational representations for $K_{2}$ and applications}
In this section we apply the concept of the right-Riemann-Liouville fractional derivative to obtain operational images for hypergeometric functions $K_{2}$. Indeed, in this section we build up the right Riemann-Liouville fractional derivative operator, which plays the role of augmenting parameters in the hypergeometric functions involved ( see \cite{sriv-karl},Chapters 4 and 5]).
In particular , we will deal with operational definitions ruled by the operators ~$\hat{D}_{x}$~and~$\hat{D}_{x}^{-1}$
~~ where ~$\hat{D}_{x}$~ denotes the derivative operator and ~$\hat{D}_{x}^{-1}$ defines the inverse of the derivative and once acting on unity yields
$$\hat{D}_{x}^{-m}1 ~=~\frac{x^{m}}{m!} , m \in N \cup \{ 0\}.$$
It is evident that $D_{x}^{-1}$~ is essentially an integral operator and the lower limit has assumed to be zero.
The following two formulas are well-known consequences of the derivative operator ~$\hat{D}_{x}$~ and the integral operator ~$\hat{D}_{x}^{-1}$~(see \cite{m-ross})
$$\hat{D}_{x}^{m}x^{\lambda} ~=~\frac{\Gamma(\lambda+1)}{\Gamma(\lambda-m+1)}x^{\lambda-m}, \eqno(3.1)$$
$$\hat{D}_{x}^{-m}x^{\lambda} ~=~\frac{\Gamma(\lambda+1)}{\Gamma(\lambda+m+1)}x^{\lambda+m}, \eqno(3.2)$$
$m \in N \cup \{ 0\} , \lambda \in C/ \{-1,-2,....\}$~.
Based on the operational relations (3.1) and (3.2) we have
$$\hat{D}_{t}^{m}\hat{D}_{u}^{-m}\{t^{\beta+m-1} u^{\gamma-
1}\}=\frac{(\beta)_{m}}{(\gamma)_{m}}\{t^{\beta-1} u^{\gamma+m-1}\}.\eqno(3.3)$$
\\
In view of the defintion of the Gaussian hypergeometric function ${}_{2}F_{1}$\cite{Sr-Ka-85}:
\\
$${}_{2}F_{1}\left[\alpha,\beta;\gamma;x\right]=\sum_{n=0}^{\infty}\frac{(\alpha)_{n}(\beta)_{n}}{(\gamma)_{n}n!}x^{n}. $$
Bin-Saad \cite{Bin-saad} introduced the following operational representations for Gaussian hypergeometric  ${}_{2}F_{1}$:
\\
{\bf{Lemma 1.}}~Let $Re(\beta)>0 ~\mbox{and}~ Re(\gamma)>0$ , then
$$\left(1-x\hat{D}_{t}u^{-1}\hat{D}_{u}^{-1}t\right)^{-\alpha}\{t^{\beta-1} u^{\gamma-1}\}
=\{t^{\beta-1} u^{\gamma-1}\}{}_{2}F_{1}\left[
\alpha,\beta ;
\gamma ; x
\right].~\eqno(3.4)$$
$$\left(1-\hat{D}_{t}\hat{D}_{x}^{-1}t\right)^{-\alpha}\{t^{\beta-1} x^{\gamma-1}\}
=\{t^{\beta-1} x^{\gamma-1}\}
{}_{2}F_{1}\left[
\alpha,\beta ;
\gamma ; x
\right].\eqno(3.5)$$
For the quadrouble hypergeometric function $K_{2}$ we introduce the following operational representations\\
{\bf{Theorem 3.1}}~Let $Re(a)>0,Re(b)>0,Re(c)>0,Re(e_{1})>0,Re(e_{2})>0,\\
Re(e_{3})>0 ~\mbox{and}~ Re(e_{4})>0$ , then
$$\left(1-x\hat{D}_{t_1}{t_2}^{-1}\hat{D}_{t_2}^{-1}t_1-y\hat{D}_{t_{1}}{t_{3}}^{-1}
\hat{D}_{t_{3}}^{-1}t_1-u\hat{D}_{t_{1}}{t_5}^{-1}\hat{D}_{t_{5}}^{-1}t_{1} \right)^{-b}$$
$$\times \left(1-z\hat{D}_{t_{1}}{t_{4}}^{-1}\hat{D}_{t_{4}}^{-1}t_{1} \right)^{-c}
\left\{ t_{1}^{a-1}t_{2}^{e_{1}-1}t_{3}^{e_{2}-1}t_{3}^{e_{3}-1}t_{5}^{e_{4}-1} \right\}$$
$$=\left\{ t_{1}^{a-1}t_{2}^{e_{1}-1}t_{3}^{e_{2}-1}t_{4}^{e_{3}-1}t_{5}^{e_{4}-1} \right\}
\times K_{2}\left[ a,a,a,a;b,b,c,b;e_{1},e_{2},e_{3},e_{4};x,y,z,u \right], \eqno(3.7)$$
$$\left(1-x\hat{D}_{t_{1}}{t_{3}}^{-1}\hat{D}_{t_{3}}^{-1}t_{1}-y\hat{D}_{t_{1}}{t_{4}}^{-1}\hat{D}_{t_{4}}^{-1}t_{1}-z\hat{D}_{t_{2}}{t_{5}}^{-1}\hat{D}_{t_{5}}^{-1}t_{2}-u\hat{D}_{t_{1}}{t_{6}}^{-1}\hat{D}_{t_{6}}^{-1}t_{1} \right)^{-a}$$
$$\times \left\{ t_{1}^{b-1}t_{2}^{c-1}t_{3}^{e_{1}-1}t_{4}^{e_{2}-1}t_{5}^{e_{3}-1}t_{6}^{e_{4}-1} \right\}$$
$$=\left\{ t_{1}^{b-1}t_{2}^{c-1}t_{3}^{e_{1}-1}t_{4}^{e_{2}-1}t_{5}^{e_{3}-1}t_{6}^{e_{4}-1} \right\}
\times K_{2}\left[ a,a,a,a;b,b,c,b;e_{1},e_{2},e_{3},e_{4};x,y,z,u \right]. \eqno(3.8)$$
{\textit{Proof.}}
 Denote , for convenience  , the left-hand side of assertion (3.7) by $I$.
 Then as a consequence of the binomial theorem , it is easily seen that:
$$I=\sum_{m,n,p,q=~0}^{\infty} \frac{(b)_{m+n+q}(c)_{p}x^{m}y^{n}z^{p}u^{q}t_{2}^{-m}t_{3}^{-n}t_{4}^{-p}t_{5}^{-q}}{m!n!p!q! }\hat{D}_{t_{1}}^{m+n+p+q}\hat{D}_{t_{2}}^{-m}\hat{D}_{t_{3}}^{-n}\hat{D}_{t_{4}}^{-p}\hat{D}_{t_{5}}^{-q}$$
$$\{t_{1}^{a+m+n+p+q-1}t_{2}^{e_{1}-1}t_{3}^{e_{2}-1}t_{4}^{e_{3}-1}t_{5}^{e_{4}-1} \}.\eqno(3.9)$$
Upon using (3.3) and considering the definition (1.1), we are led finally to right-hand side the assertion (3.7). The proof of the assertion (3.8) runs parallel to that of the assertion (3.7) then we skip the details.
\\
\\
One of the advantages offered by the use of the operational images (3.7) and (3.8) is the possibility of establishing finite sums for the quadruple hypergeometric series $K_{2}$.
 First,let, for convenience,\\
$\hat{M}=x\hat{D}_{t_{1}}{t_{3}}^{-1}\hat{D}_{t_{3}}^{-1}t_{1}-y\hat{D}_{t_{1}}{t_{4}}^{-1}\hat{D}_{t_{4}}^{-1}t_{1}-z\hat{D}_{t_{2}}{t_{5}}^{-1}\hat{D}_{t_{5}}^{-1}t_{2}-u\hat{D}_{t_{1}}{t_{6}}^{-1}\hat{D}_{t_{6}}^{-1}t_{1}$
\\
$\hat{N}=x\hat{D}_{t_{1}}{t_{2}}^{-1}\hat{D}_{t_{2}}^{-1}t_{1}-y\hat{D}_{t_{1}}{t_{3}}^{-1}\hat{D}_{t_{3}}^{-1}t_{1}-u\hat{D}_{t_{1}}{t_{5}}^{-1}\hat{D}_{t_{5}}^{-1}t_{1}$
\\
$\hat{R}=z\hat{D}_{t_{1}}{t_{4}}^{-1}\hat{D}_{t_{4}}^{-1}t_{1}$
\\
By employing (3.8) and taking advantage from
the identity
$$\left[1-\left(1-\hat{M}\right)\right]^{n}\left\{ t_{1}^{b-1}t_{2}^{c-1}t_{3}^{e_{1}-1}t_{4}^{e_{2}-1}t_{5}^{e_{3}-1}t_{6}^{e_{4}-1} \right\}=\left[\hat{M}\right]^{n}\left\{ t_{1}^{b-1}t_{2}^{c-1}t_{3}^{e_{1}-1}t_{4}^{e_{2}-1}t_{5}^{e_{3}-1}t_{6}^{e_{4}-1} \right\}$$
and the definition of the general triple hypergeometric series $F^{(3)}[x,y,z]$(see \cite{Sr-Ka-85}, we can state that
$$\sum_{r=0}^{n}(-1)^{r}\left(\begin{array}{l}
n \\
r \\
\end{array}\right)
K_{2}\left[ -r,-r,-r,-r;b,b,c,b;e_{1},e_{2},e_{3},e_{4};x,y,z,u \right]$$
$$=\frac{(b)_{n}}{(e_{1})_{n}} F^{(3)}
\left[\begin {array}{l}
-n~ , 1-e_{1}-n~::--~;c~;--;  \\
                                    \frac{y}{x}, \frac{z}{x}, \frac{u}{x}\\
-------::e_{2}~;e_{3}~;e_{4}~;\\
\end{array}\right]~.\eqno(3.10)$$\\
Similarly,in view of (3.7),the definition of  Appell hypergeometric series $F_{4}[x,y]$(see [22,p.23(5)]) and the use of the following identity\\
$$\left[1-\left(1-\hat{N}\right)\right]^{n}\left[1-\left(1-\hat{R}\right)\right]^{m}\left\{ t_{1}^{a-1}t_{2}^{e_{1}-1}t_{3}^{e_{2}-1}t_{4}^{e_{3}-1}t_{5}^{e_{4}-1} \right\}$$
$$=\left[\hat{N}\right]^{n}\left[\hat{R}\right]^{m}\left\{ t_{1}^{a-1}t_{2}^{e_{1}-1}t_{3}^{e_{2}-1}t_{4}^{e_{3}-1}t_{5}^{e_{4}-1} \right\}$$
\\
allows to conclude that
$$\sum_{r=0}^{n}\sum_{s=0}^{m}(-1)^{r+s}
\left(\begin{array}{l}
n \\
r \\
\end{array}\right)
\left (\begin{array}{l}
m \\
s \\
\end{array}\right)
K_{2}\left[a,a,a,a, -r,-r,-s,-r;e_{1},e_{2},e_{3},e_{4};x,y,z,u \right]$$
$$=\frac{(a)_{m+n}}{(e_{1})_{n}(e_{3})_{m}} F_{4}\left[-n,1-e_{1}-n;e_{2},e_{4};\frac{y}{x},\frac{u}{x}\right]\eqno(3.11)$$
Also, by virtue of the identity
\\
$$\left[\left(1-\hat{N}\right)\right]^{-b}\left[\left(1-\hat{R}\right)\right]^{-b}\left\{ t_{1}^{a-1}t_{2}^{e_{1}-1}t_{3}^{e_{2}-1}t_{4}^{e_{3}-1}t_{5}^{e_{4}-1} \right\}$$
$$=\left[1-\hat{N}-\hat{R}+\hat{N}\hat{R}\right]^{-b}\left\{ t_{1}^{a-1}t_{2}^{e_{1}-1}t_{3}^{e_{2}-1}t_{4}^{e_{3}-1}t_{5}^{e_{4}-1} \right\}$$
\\
we find the following decomposition formula:\\
$$\sum_{s,k,r=0}^{\infty} \frac{(a)_{2s+2k+2r} (b)_{s+k+r}}{(e_{1})_{s}(e_{2})_{k}(e_{3})_{s+k+r}(e_{4})_{r}} (-xz)^{s}(-yz)^{k}(-uz)^{r}$$
$$\times F_{C}^{(4)}\left[a+2s+2k+2r,b+s+k+r;e_{1},e_{2},e_{3},e_{4};x,y,z,u\right]$$
$$K_{2}\left[ a,a,a,a;b,b,b,b;e_{1},e_{2},e_{3},e_{4};x,y,z,u \right],\eqno(3.12)$$\\
where $F_{C}^{(4)}$ is Lauricella function of four variables[22,p.33(3)].\\
\\
\\
Moreover, in view of the composition rule\\
$\left[\left(1-\hat{M}\right)\right]^{-2a}=\left[\left(1-\hat{M}\right)\right]^{-a}\left[\left(1-\hat{M}\right)\right]^{-a}=\left[1-2\hat{M}+\hat{M}^{2}\right]^{-a}$\\
and the operational image (3.8), we easily find for the  function $K_{2}$ that\\
$$K_{2}\left[2a,2a,2a,2a,b,b,c,b;e_{1},e_{2},e_{3},e_{4};x,y,z,u\right]$$
$$=\sum_{m,n,p,q,s,r,k,l,h,t=0}^{\infty}\frac{(a)_{M}(b)_{N}(c)_{P}}{(e_{1})_{Q^{(1)}}(e_{2})_{Q^{(2)}}(e_{3})_{Q^{(3)}}(e_{4})_{Q^{(4)}}}x^{2m}y^{2n}z^{2p}u^{2q}(2xy)^{s}(2xz)^{k}(2xu)^{r}$$
$$(2yz)^{l}(2yu)^{h}(2zu)^{t} \times K_{2}\left[a+M,a+M,a+M,a+M,b+N,b+N,c+N,b+N;\right.$$
$$\left. e_{1}+Q^{(1)},e_{2}+Q^{(2)},e_{3}+Q^{(3)},e_{4}+Q^{(4)};x,y,z,u \right].\eqno(3.13)$$
where
$$\begin{pmatrix}
M:=m+n+p+q+s+k+r+l+h+t,P:=2p+k+l+t,&\\
N:=2(m+n+q+s+r+h)+k+l+t,\\
Q^{(1)}:=2m+s+k+r,&\\
Q^{(2)}:=2n+s+l+h,Q^{(3)}:=2p+k+l+t,Q^{(4)}:=2q+r+h+t&
\end{pmatrix}.\eqno(3.14)$$

Finally, let us stress that the schema suggested in sections 2 and 3 can be applied to find linear independent solutions , symbolic operational images , finite sums and decomposition formulas for other quadruple hypergeometric functions.

\end{document}